\documentclass[10pt,conference,letterpaper,twocolumn]{ieeeconf} 
\usepackage{graphicx} 

\IEEEoverridecommandlockouts
\usepackage{comment}

\usepackage
{
    graphicx,
    amssymb,
    amsmath,
    dsfont, 
    xcolor,
    mathtools,
    authblk,
    nicefrac,
    array,
    enumerate,
    siunitx,
    tikz,
}
\usepackage{algorithm}
\usetikzlibrary{shapes,arrows}
\usetikzlibrary{calc} 

\usepackage[utf8]{inputenc}

\usepackage[pdffitwindow=false,
            plainpages=false,
            pdfpagelabels=true,
            pdfpagemode=UseOutlines,
            pdfpagelayout=SinglePage,
            bookmarks=false,
            colorlinks=true,
            hyperfootnotes=false,
            linkcolor=blue,
            urlcolor=blue!30!black,
            citecolor=green!50!black]{hyperref}

\usepackage{caption}

\graphicspath{{resources/}}

\usepackage{cleveref}
\makeatletter
\newcommand\RedeclareMathOperator{%
  \@ifstar{\def\rmo@s{m}\rmo@redeclare}{\def\rmo@s{o}\rmo@redeclare}%
}
\newcommand\rmo@redeclare[2]{%
  \begingroup \escapechar\m@ne\xdef\@gtempa{{\string#1}}\endgroup
  \expandafter\@ifundefined\@gtempa
     {\@latex@error{\noexpand#1undefined}\@ehc}%
     \relax
  \expandafter\rmo@declmathop\rmo@s{#1}{#2}}
\newcommand\rmo@declmathop[3]{%
  \DeclareRobustCommand{#2}{\qopname\newmcodes@#1{#3}}%
}
\@onlypreamble\RedeclareMathOperator
\makeatother

\newcommand{\N}{\mathds{N}}
\newcommand{\R}{\mathds{R}}
\newcommand{\Rp}{\R_{\geq0}}

\newcommand{\Impl}{\Longrightarrow}

\newcommand{\fa}{\ \forall \, }
\newcommand{\ex}{\ \exists \, }

\newcommand{\rbl}{\left (}
\newcommand{\rbr}{\right )}

\newcommand{\al}{\left \langle}
\newcommand{\ar}{\right \rangle}
\newcommand{\bl}{\left |}
\newcommand{\br}{\right |}
\newcommand{\nl}{\left\|}
\newcommand{\nr}{\right\|}
\newcommand{\cbl}{\left\lbrace }
\newcommand{\cbr}{\right\rbrace }
\newcommand{\Abs}[1]{\bl #1 \br}
\newcommand{\Norm}[2][ ]{\nl #2 \nr_{#1}}
\newcommand{\SNorm}[1]{\Norm[\infty]{#1}}
\newcommand{\LNorm}[2][2]{\Norm[L^{#1}]{#2}}

\newcommand{\setdef}[2]{\cbl\ #1\ \left|\ \vphantom{#1} #2\ \right.\cbr}

\newcommand{\GL}{\text{GL}}
\newcommand{\cB}{\mathcal{B}}
\newcommand{\cD}{\mathcal{D}}

\newcommand{\cU}{\mathcal{U}}

\newcommand{\cF}{\mathcal{F}}
\newcommand{\cG}{\mathcal{G}}
\newcommand{\cP}{\mathcal{P}}
\newcommand{\cC}{\mathcal{C}}
\newcommand{\cN}{\mathcal{N}}
\newcommand{\cY}{\mathcal{Y}}

\newcommand{\cT}{\mathcal{T}}

\DeclareMathOperator*{\rf}{ref}
\DeclareMathOperator*{\esssup}{ess\,sup}

\DeclareMathOperator*{\loc}{loc}

\newcommand{\me}{\mathrm{e}}

\newcommand{\con}{\mathcal{C}}

\newcommand{\oT}{\mathbf{T}}

\RedeclareMathOperator*{\Im}{Im}
\RedeclareMathOperator*{\Re}{Re}
\renewcommand{\phi}{\varphi}

\renewcommand{\d}{\ \text{d}}

\newcommand{\dd}[2][ ]{\tfrac{\text{\normalfont d}#1}{\text{\normalfont d}#2}}

\makeatletter
\@ifundefined{ifCustomTheorems}{}
{
    \setcounter{section}{0}
    \counterwithout{equation}{section}
    \counterwithout{figure}{section}
    \newtheorem{definition}{Definition}[section]
    \theoremstyle{definition}
    \newtheorem{remark}[definition]{Remark}\Crefname{remark}{Remark}{Remarks}
    \Crefname{algo}{Algorithm}{Algorithms}
    \Crefname{example}{Example}{Examples}
    \theoremstyle{plain}
    \Crefname{prop}{Proposition}{Propositions}
    \newtheorem{corollary}[definition]{Corollary}\Crefname{corollary}{Corollary}{Corollaries}
    \Crefname{assertion}{Assertion}{Assertions}
    \newtheorem{theorem}[definition]{Theorem}\Crefname{theorem}{Theorem}{Theorems}
    \newtheorem{lemma}[definition]{Lemma}\Crefname{lemma}{Lemma}{Lemmata}
}
\makeatother

{\left(\begin{smallmatrix}}%
{\end{smallmatrix}\right)}%

\newenvironment{smallbmatrix}%
{\left[\begin{smallmatrix}}%
{\end{smallmatrix}\right]}%

\newtheorem{theorem}{Theorem}[section]

\newtheorem{definition}{Definition}[section]
\newtheorem{remark}{Remark}[section]
\newtheorem{lemma}{Lemma}[section]
\newtheorem{corollary}{Corollary}[section]
\Crefname{definition}{Definition}{Definitions}
\Crefname{assumption}{Assumption}{Assumptions} 
\Crefname{lemma}{Lemma}{Lemmata}
\Crefname{corollary}{Corollary}{Corollaries}
\Crefname{remark}{Remark}{Remarks}
\Crefname{theorem}{Theorem}{Theorems}

\title{
On model predictive control with sampled-data input for output tracking with prescribed performance
}
\author{Dario Dennstädt{$^*$}{$^\dagger$} \and Lukas Lanza{$^\dagger$} \and Karl Worthmann{$^\dagger$}
\thanks{We gratefully acknowledge funding by the Deutsche Forschungsgemeinschaft (DFG, German Research Foundation) – Project-IDs 471539468 and 507037103; Lukas Lanza is also grateful for the support from the Carl Zeiss Foundation (VerneDCt -- Project No.\ 2011640173).}
\thanks{$^{\dagger}$ Technische Universität Ilmenau, Institute of Mathematics,  Optimization-based Control Group, Weimarer Stra{\ss}e~25, 98693~Ilmenau, Germany
\texttt{ \{lukas.lanza, karl.worthmann\}@tu-ilmenau.de}}
\thanks{$^{*}$ Universit\"at Paderborn, Institut f\"ur Mathematik, Warburger Stra{\ss}e~100, 33098~Paderborn, Germany (\texttt{dario.dennstaedt@uni-paderborn.de})} }

\begin{document}

\maketitle
\begin{abstract}
We propose a model predictive control (MPC)
scheme with sampled-data input
which ensures output-reference tracking within prescribed error bounds for relative-degree-one systems. 
Hereby, we explicitly deduce bounds 
on the required maximal control input
and sampling frequency such that the MPC scheme is both initially and recursively feasible.
A key feature of the proposed approach is that neither terminal conditions nor a sufficiently-large prediction horizon are imposed,
rendering the MPC scheme computationally efficient. 
We illustrate the MPC algorithm via a numerical example of a torsional oscillator.
\end{abstract}

\section{Introduction}
Model predictive control (MPC) has gained widespread recognition due to its ability to effectively deal with nonlinear multi-input multi-output systems while adhering to control and state constraints,
see the textbooks~\cite{grune2017nonlinear,rawlings2017model} and the references therein.
Since MPC is based  on the iterative solution of optimal control problems (OCPs) over a finite prediction horizon,
guaranteeing their initial and recursive feasibility is key to  ensure proper functioning of MPC.
W.r.t.\ recursive feasibility, often either terminal conditions~\cite{chen1998quasi} or a combination of cost controllability~\cite{CoroGrun20} and a sufficiently long prediction horizon are used, 
see, e.g., \cite{EsteWort20} and the references therein for discrete-time systems. 
The use of terminal constraints may increase the computational effort and substantially reduce the domain of attraction, see, e.g.,~\cite{gonzalez2009enlarging}. 
This becomes even more involved in the presence of time-varying state or output constraints~\cite{manrique2014mpc}.

Recently, a novel MPC scheme was proposed in~\cite{berger2019learningbased} to overcome these restrictions by
invoking structural properties of the system class in consideration to show both initial and recursive feasibility. 
This approach was further developed in \cite{BergDenn21,BergDenn23b} and, using
the combination of a high-gain property based on a well-defined relative degree and input-to-state stable internal dynamics,
it allows for output reference tracking within prescribed time-varying error bounds of continuous-time systems.
This also distinguishes it from the output regulation problem, on which other MPC approaches focus, see, e.g.,~\cite{limon2018nonlinear}, 
where discrete-time systems are considered.
The output regulation regulation problem was considered in the presence of time-invariant constraints in~\cite{kohler2021constrained},
where suitable stabilizability and detectability conditions and a sufficiently long prediction horizon are used to ensure constraint satisfaction.

The underlying idea of the approach from~\cite{berger2019learningbased,BergDenn21,BergDenn23b} is closely intertwined with the adaptive, high-gain funnel controller~\cite{IlchRyan02b}, see, e.g., the recent survey paper~\cite{berger2023funnel}, which also guarantees output tracking within given bounds by continuously adapting the applied control signal based on continuously available measurements.
However, both in practical applications as well as in simulations, system measurements and control input signals are 
typically only given at discrete time instances. 
Consequently, the assumption underlying funnel control, namely that the signal is continuously available, is not met.
In the recent work~\cite{LanzaDenn23sampled} this shortcoming was addressed. It was shown that 
the control objective of ensuring output tracking with predefined error boundaries can be achieved  
by a sampled-data feedback controller,
that receives system measurements only at uniformly sampled discrete time instances and yields a piecewise constant control signal.
The latter is usually called \emph{sampled-data control}.

In~\cite{yuz2005sampled}, it is shown how to obtain a sampled-data model approximation for continuous-time systems, where the mismatch between the solutions scales with the sampling time.
Based on the concept of control barrier functions (CBFs) for continuous-time systems, in~\cite{breeden2021control} sampled-data CBFs are utilized to ensure safe sets to be forward invariant.
In, e.g.,~\cite{geromel2021sampled,worthmann2015unconstrained} a sampled-data MPC scheme for continuous-time systems is developed, respectively. 

In this paper, we show that it is possible to achieve output tracking with prescribed performance using MPC with sampled-data control.
Although simulations suggest that the MPC scheme from~\cite{berger2019learningbased,BergDenn21} can be implemented using sampled-data inputs,
theoretical results, so far, are missing.
We rigorously prove that the novel MPC algorithm proposed in the present paper is initially and recursively feasible.
Furthermore, we derive explicit bounds on the maximal required control input as well as on the sufficiently small step length. 
To this end, we combine results from~\cite{LanzaDenn23sampled} on sampled-data adaptive feedback control and from the continuous-time MPC case~\cite{BergDenn21}.

The article is organized as follows.
In \Cref{Sec:SystemClassAndControlObjective}, we introduce the system class 
and define the 
control objective.
In \Cref{Sec:MPCAlgorithm} we propose an MPC algorithm, which achieves the control objective by using sampled-data inputs only.
Its feasibility is proven in our main result \Cref{Th:FunnelMPC}.
We illustrate the MPC algorithm via a simulation in \Cref{Sec:Simulation}. 
Conclusions and an outlook are given in \Cref{Sec:Conclusion}.

\textbf{Nomenclature.}
$\N$, $\R$ denote natural and real numbers, respectively.
$\N_0:=\N\cup\{0\}$ and $\Rp:=[0,\infty)$.
$\Norm{x}:=\sqrt{\al x,x\ar}$~denotes the Euclidean norm of $x\in\R^n$, and $\cB_v :=\setdef{ x \in \R^n}{\|x\| \le v}$.
$\Norm{A}$ denotes the induced operator norm
$\Norm{A}:=\sup_{\Norm{x} = 1}\Norm{Ax}$ for $A\in\R^{n\times n}$.
$\GL_n(\R)$ is the group of invertible $\R^{n\times n}$ matrices.
$\con^p(V,\R^n)$ is the linear space of $p$-times continuously  differentiable
functions $f:V\to\R^n$, where $V\subset\R^m$ and $p\in\N_0\cup \{\infty\}$.
$\con(V,\R^n):=\con^0(V,\R^n)$.
On an interval $I\subset\R$,  $L^\infty(I,\R^n)$ denotes the space of measurable and essentially bounded
functions $f: I\to\R^n$ with norm $\SNorm{f}:=\esssup_{t\in I}\Norm{f(t)}$,
$L^\infty_{\text{loc}}(I,\R^n)$ the set of measurable and locally essentially bounded functions, and $L^p(I,\R^n)$
the space of measurable and $p$-integrable functions with norm $\LNorm[p]{\cdot}$ and with $p\ge 1$.
Furthermore, $W^{k,\infty}(I,\R^n)$ is the Sobolev space of all $k$-times weakly differentiable functions
$f:I\to\R^n$ such that $f,\dots, f^{(k)}\in L^{\infty}(I,\R^n)$.

\section{System class and control objective} \label{Sec:SystemClassAndControlObjective}
In this section, we introduce the class of systems to be controlled and define the control objective precisely.
\subsection{System class} \label{Ssec:SysClass}
We consider nonlinear multi-input multi-output systems
\begin{equation} \label{eq:Sys}
    \begin{aligned}
    \dot{y}(t) &= f \big(\oT(y)(t) \big) + g(\oT(y)(t))u(t) \\
    y|_{[t_0-\sigma,t_0]} &= y^0  \in \con([t_0-\sigma,t_0],\R^m),
    \end{aligned}
\end{equation}
with initial time $t_0\ge 0$, ``memory'' $\sigma \ge 0$, initial trajectory $y^0$,
control input~$u\in L^\infty_{\loc}(\Rp, \R^m)$, 
and output $y(t)\in\R^m$ at time $t\geq t_0$.
Note that $u$ and $y$ have the same dimension~${m\in\N}$.
The system consists of the  nonlinear functions 
$f \in\con(\R^q,\R^m)$, $g \in \con(\R^q , \R^{m \times m})$,
and the  nonlinear operator $\oT:\con([-\sigma,\infty),\R^m)\to L^\infty_{\loc}(\Rp,\R^q)$.
The operator~$\oT$ is causal, locally Lipschitz, and satisfies a bounded-input bounded-output property. 
It is characterized in detail in the following definition.
\begin{definition} \label{Def:OperatorClass} 
For $m,q\in\N$ and $\sigma\geq 0$, the set $\cT_\sigma^{m,q}$ denotes the class of operators $\textbf{T}:
\con([t_0-\sigma,\infty),\R^m) \to L^\infty_{\loc} (\Rp, \R^{q})$
for which the following properties hold:
\begin{itemize}
    \item\textit{Causality}:  $\fa y_1,y_2\in\con([t_0-\sigma,\infty),\R^m)$  $\fa t\geq t_0$:
    \[
        y_1\vert_{[t_0-\sigma,t]} = y_2\vert_{[t_0-\sigma,t]}
        \ \Impl\ 
        \textbf{T}(y_1)\vert_{[t_0,t]}=\textbf{T}(y_2)\vert_{[t_0,t]}.
    \]
    \item\textit{Local Lipschitz}: 
    $\fa t \ge t_0 $ $\fa y \in \con([t_0-\sigma,t] , \R^m)$ 
    $\ex \Delta, \delta, c > 0$ 
    $\fa y_1, y_2 \in \con([t_0-\sigma,\infty) , \R^m)$ with
    $y_1|_{[t_0-\sigma,t]} = y_2|_{[t_0-\sigma,t]} = y $ 
    and $\Norm{y_1(s) - y(t)} < \delta$,  $\Norm{y_2(s) - y(t)} < \delta $ for all $s \in [t,t+\Delta]$:
    \[
     \hspace*{-2mm}   \esssup_{\mathclap{s \in [t,t+\Delta]}}  \Norm{\textbf{T}(y_1)(s) \!-\! \textbf{T}(y_2)(s) }  
        \!\le\! c \ \sup_{\mathclap{s \in [t,t+\Delta]}}\ \Norm{y_1(s)\!-\! y_2(s)}\!.
    \] 
    \item\textit{Bounded-input bounded-output (BIBO)}:
    $\fa c_0 > 0$ $\ex c_1>0$  $\fa y \in \con([t_0-\sigma,\infty), \R^m)$:
    \[
    \sup_{t \in [t_0-\sigma,\infty)} \Norm{y(t)} \le c_0 \ 
    \Impl \ \sup_{t \in [t_0,\infty)} \Norm{\textbf{T}(y)(t)}  \le c_1.
    \]
\end{itemize}
\end{definition}
Since the operator~$\oT$ acts on the whole output trajectory, the causality property of~\Cref{Def:OperatorClass}
enforces that the system does not depend on future system states. 
The second condition (locally Lipschitz)
is more of a technical nature to guarantee existence and uniqueness of solutions.
Finally, the BIBO property ensures that the system remains bounded as long as the system output does.
Note that using the operator~$\oT$ many physical phenomena such as \emph{backlash},
and \emph{relay hysteresis}, and \emph{nonlinear time delays}
can be modeled, where $\sigma \ge 0$ corresponds to the initial delay, cf.~\cite[Sec.~1.2]{BergIlch21}.
\begin{remark} \label{Rem:StateRepresentation}
Consider a nonlinear control affine system
\begin{equation}\label{eq:SysWithState}
\begin{aligned}
    \dot{x}(t)  & = F(x(t)) + G(x(t)) u(t),\quad x(t_0)=x^0\in\R^n,\\
    y(t)        & = H(x(t)),
\end{aligned} 
\end{equation}
with $t_0\in\Rp$, $x^0\in\R^n$, and nonlinear functions $F\in\con^1(\R^n,\R^n)$, $G \in \con^1(\R^n,\R^{n \times m})$, and  $H\in\con^2(\R^n,\R^m)$.
Under assumptions provided in~\cite[Cor.~5.6]{ByrnIsid91a}, there exists a diffeomorphism $\Phi:\R^n\to\R^n$ which induces a coordinate transformation 
putting the system~\eqref{eq:SysWithState} into the form~\eqref{eq:Sys} 
with new coordinates~$(y,\eta) \in \R^m \times \R^{n-m}$ (output and internal state)
for appropriate functions $f$, $g$, operator~$\oT$, and $\sigma=0$.
In this case $\oT$ is the solution operator of the internal dynamics of the transformed system.
As in~\cite{BergDenn21}, exact knowledge about the coordinate transformation and computation of the diffeomorphism~$\Phi$ is not required to apply \Cref{Algo:DiscreteFMPC} to the system~\eqref{eq:SysWithState} -- 
merely the existence of~$\Phi$ has to be assumed as a mean for the proofs.
\end{remark}

Invoking the requirements for the operator in~\Cref{Def:OperatorClass}, we formally introduce the system class. 
\begin{definition} \label{Def:system-class}
    For $m \in \N$, a system~\eqref{eq:Sys} belongs to the system class $\cN^{m}$, written ${(f,g,\oT)
    \in\cN^{m}}$, if, for some $q\in\N$ and $\sigma \geq0$, the following holds:
    ${f\in\con(\R^q ,\R^m})$, ${\oT\in\cT^{m,q}_{\sigma}}$, and 
    the function~$g$  is strictly positive definite, that is 
    for all $\xi \in \R^{q}$  and for all $z \in \R^m \setminus \{0\} $ 
    \begin{equation} \label{Ass:g_sign_definite}
         \al z, g({\xi}) z \ar > 0.
    \end{equation}
\end{definition}

For $t\geq0 $ and a control function $u\in L_{\loc}^\infty(\Rp,\R^m)$, the system~\eqref{eq:Sys} 
has a solution in the sense of \textit{Carath\'{e}odory}, meaning there exists a function $y:[t_0-\sigma,\omega)\to\R^m$, $\omega>t_0$, 
with $y|_{[t_0-\sigma,t_0]} = y^0  \in \con([t_0-\sigma,t_0],\R^m)$
and $y|_{[t_0,\omega)}$ is absolutely continuous and satisfies the ODE in~\eqref{eq:Sys} for almost all $t\in[t^0,\omega)$.
A solution $y$ is called \textit{maximal}, if it has no right extension that is also a solution.
A maximal solution is called \textit{response} associated with $u$ and denoted by~$y(\cdot;t_0,y^0,u)$.
Note that in the case~$\sigma=0$, we mean by $y|_{[t_0-\sigma,t_0]}$ the evaluation of the function at $t_0$,
i.e., $y|_{[t_0-\sigma,t_0]}=y(t_0)$, and refer to the vector space $\R^m$ when using the notation $\con([t_0-\sigma,t_0],\R^m)$.

\subsection{Control objective}\label{Ssec:ContrObj}
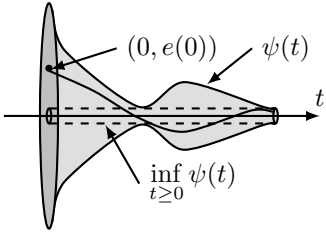
\begin{figure}[h]
  \begin{center}
\begin{tikzpicture}[scale=0.30]
\tikzset{>=latex}
  \filldraw[color=gray!25] plot[smooth] coordinates {(0.15,4.7)(0.7,2.9)(4,0.4)(6,1.5)(9.5,0.4)(10,0.333)(10.01,0.331)(10.041,0.3) (10.041,-0.3)(10.01,-0.331)(10,-0.333)(9.5,-0.4)(6,-1.5)(4,-0.4)(0.7,-2.9)(0.15,-4.7)};
  \draw[thick] plot[smooth] coordinates {(0.15,4.7)(0.7,2.9)(4,0.4)(6,1.5)(9.5,0.4)(10,0.333)(10.01,0.331)(10.041,0.3)};
  \draw[thick] plot[smooth] coordinates {(10.041,-0.3)(10.01,-0.331)(10,-0.333)(9.5,-0.4)(6,-1.5)(4,-0.4)(0.7,-2.9)(0.15,-4.7)};
  \draw[thick,fill=lightgray] (0,0) ellipse (0.4 and 5);
  \draw[thick] (0,0) ellipse (0.1 and 0.333);
  \draw[thick,fill=gray!25] (10.041,0) ellipse (0.1 and 0.333);
  \draw[thick] plot[smooth] coordinates {(0,2)(2,1.1)(4,-0.1)(6,-0.7)(9,0.25)(10,0.15)};
  \draw[thick,->] (-2,0)--(12,0) node[right,above]{\normalsize$t$};
  \draw[thick,dashed](0,0.333)--(10,0.333);
  \draw[thick,dashed](0,-0.333)--(10,-0.333);
  \node [black] at (0,2) {\textbullet};
  \draw[->,thick](4,-3)node[right]{\normalsize$\inf\limits_{t \ge 0} \psi(t)$}--(2.5,-0.4);
  \draw[->,thick](3,3)node[right]{\normalsize$(0,e(0))$}--(0.07,2.07);
  \draw[->,thick](9,3)node[right]{\normalsize$\psi(t)$}--(7,1.4);
\end{tikzpicture}
\end{center}
 \vspace*{-2mm}
 \caption{Error evolution in a funnel $\mathcal F_{\psi}$ with boundary $\psi(t)$; the figure is based on~\cite[Fig.~1]{BergLe18a}, edited for present purpose.}
 \label{Fig:funnel}
 \vspace*{-2mm}
 \end{figure}

The control objective is that the output~$y$ of system~\eqref{eq:Sys} follows a given reference~$y_{\rm ref}$ with predefined accuracy.
To be more
precise, the tracking error ~$t\mapsto e(t):=y(t)-y_{\rf}(t)$ shall evolve within the prescribed
performance funnel
\begin{align*}
    \cF_\psi= \setdef{(t,e)\in \Rp\times\R^{m}}{\Norm{e} \leq \psi(t)},
\end{align*}
or formulated in a different way, the output~$y(t)$ should at every time instance $t\geq t_0$ belong to the set
\begin{align}\label{eq:Def-Dt}
    \cD_t := \setdef
                    {y\in\R^m}
                    {\Norm{y - y_{\rf}(t)} \leq \psi(t)}.
\end{align}
The performance set~$\cF_\psi$ is determined by the choice of the function~$\psi$
belonging  to
\begin{align*}
    \cG:=\setdef
        {\psi\in W^{1,\infty}(\Rp,\R)}
        {
         \inf_{s \ge 0} \psi(s) > 0
        },
\end{align*}
see also Figure~\ref{Fig:funnel}.
Note that keeping the tracking error in~$\mathcal{F}_\psi$ does not mean asymptotic convergence to zero.
Moreover, the funnel boundary is not necessarily monotonically decreasing.
The specific application usually dictates the constraints on the tracking error and thus
indicates suitable choices for~$\psi$.

\section{Sampled-data MPC} \label{Sec:MPCAlgorithm}
We aim to develop an MPC algorithm, which achieves the aforementioned control objective.
In contrast to the MPC scheme proposed in~\cite{berger2019learningbased,BergDenn21}, 
the space of admissible controls is restricted to step functions, i.e., the control signal can only change finitely often between two sampling instances. 
In other words, we use sampled-data control.
To introduce the control scheme properly, we formally define step functions in the following definition.

\begin{definition} \label{Def:PartitionAndStepFunction}
    Let $I\subset\Rp$ be an interval of the form $I=[a,b]$ with $b>a$ or $I=[a,\infty)$.
    We call a strictly increasing sequence $\cP=(t_k)_{k\in\N_0}$ with $\lim_{k\to\infty}t_k=\infty$ 
    and ${t_0 = a}$ a~\textit{partition} of $I$.
    The norm of~$\cP$ is defined as $\Abs{\cP}:= \sup\setdef{t_{i+1}-t_{i}}{i\in\N_0}$.
    A function~$f:I\to\R^m$ is called~\textit{step function with partition}~$\cP$ if $f$ is constant on every
    interval~$[t_i,t_{i+1})\cap I$ for all~$i\in\N_0$. 
    We denote the space of all step functions on~$I$ with partition~$\cP$ by~$\cT_{\cP}(I,\R^m)$.
\end{definition}

Note that in the case of finite intervals $I=[a,b]$ with $b>a$,
\Cref{Def:PartitionAndStepFunction} can also be formulated 
using finite sequences  $\cP=(t_k)_{k=0}^N$ with $N\in\N$ and $t_N=b$.
However, using infinite sequences every partition~$\cP$ of $[a,\infty)$
is also a partition of $[a,b]$ for all $b>a$.
Using this fact will simplify formulating our results.
Changing the control signal is restricted to the time instances $t_k$.
Further note that \Cref{Def:PartitionAndStepFunction} allows for non-uniform step length, i.e., for $\tau_i = |t_{i-1} - t_i|$ we allow $\tau_k \neq \tau_j$ for $k \neq j$, where $i,j,k \in \N$.
However, in practice, a uniform step length often will be used.

Before formulating an MPC algorithm which achieves the control objective described in~\Cref{Ssec:ContrObj},
we introduce the class of admissible stage-costs.
Let $\tilde \ell\in\cC(\Rp\times\R^m\times\R^m,\Rp)$.
Then, we define the stage-cost~$\ell$ piecewise by
\begin{align}\label{eq:stageCostFunnelMPC}
    \ell(t,y,u)\!=\!
    \begin{dcases}\!
        \tilde{\ell}(t,y, u),
            & \Norm{e(y,t)} \!\leq\! \psi(t), \\
        \!\infty,&\text{else},
    \end{dcases}
\end{align}
where $e(y,t):= y-y_{\rf}(t)$. 
A suitable choice is, for instance,
\begin{equation} \label{eq:QuadraticCost}
\begin{aligned}
    \ell(t,y,u)=
    \begin{cases}
        \| e(y,t)\|^2 + \lambda_u \|u\|^2,
            & \Norm{e(y,t)} \leq  \psi(t)\\
        \infty, &\text{else},
    \end{cases}
\end{aligned}
\end{equation}
where~$\lambda_u \ge 0$ is a design parameter.
Note that in contrast to the MPC scheme investigated in~\cite{berger2019learningbased,BergDenn21},
we allow for a fairly large class of cost functions since the function~$\tilde{\ell}$ can be freely chosen by the user. Moreover,
the stage-cost~\eqref{eq:stageCostFunnelMPC} allows the error to be ``on the funnel boundary'', i.e., we allow for $\|y(t) - y_{\rm ref}(t) \| = \psi(t)$.

Invoking stage-costs like in~\eqref{eq:stageCostFunnelMPC}, we propose the
sampled-data MPC \Cref{Algo:DiscreteFMPC},
where the input is restricted to piecewise constant step functions with given step length.
\begin{algorithm}[H]\caption{Sampled-data MPC}\label{Algo:DiscreteFMPC}\ \\
    \textbf{Given:} System~\eqref{eq:Sys}, reference~$y_{\rf}\in
    W^{1,\infty}(\Rp,\R^{m})$, funnel function $\psi\in\cG$, control bound $u_{\max}>0$, maximal step length~$\tau>0$ of the control signal, and
    initial data $y^0\in\con([t_0-\sigma,t_0],\R^m)$.\\
    \textbf{Set} the time shift $\delta >0$, the prediction horizon $T\geq\delta$, 
    initialize the current time~$\hat{t}=t_0$, and
    choose a partition $\cP=(t_k)_{k\in\N_0}$ of the interval~$[t_0,\infty)$ with~$\Abs{\cP}\leq\tau$ 
    and which contains $(t_0 + k\delta)_{k\in\N_0}$ as a subsequence.
    \\
    \textbf{Steps:}
    \begin{enumerate}[1.]
    \item\label{agostep:FMPCFirst} Obtain a measurement of the output $y$ of~\eqref{eq:Sys} 
    on the interval~$[t_0-\sigma,\hat t]$ and set $\hat y := y|_{[t_0-\sigma,\hat{t}]}$.
    \item Compute a solution $u^{\star}\in \cT_{\cP}([t_k,t_k+T],\R^{m})$ of
    \begin{equation}\label{eq:DiscreteFMPCOCP}
            \mathop
            {\operatorname{minimize}}_{\substack{u\in \cT_{\cP}([t_k,t_k+T],\R^{m}),\\\SNorm{u}\leq u_{\max}}}  \quad
            \int_{t_k}^{t_k+T} \ell(t,y(t;\hat{t},\hat{y},u),u)
            \d t .
    \end{equation}
    \item\label{algo:item:BoundedInput}  Apply the sampled-data 
    control signal 
    $\mu:[\hat{t},\hat{t}+\delta)\times\con([t_0-\sigma,\hat{t}],\R^m)\to\R^m$, defined by 
    \begin{equation}\label{eq:FMPC-fb}
        \mu(t,\hat y) =u^{\star}(t)
    \end{equation}
        to system~\eqref{eq:Sys}.
        Increase~$\hat{t}$ by~$\delta$ and go to Step~\ref{agostep:FMPCFirst}.
    \end{enumerate}
\end{algorithm}

\begin{remark}
        If the system is given as nonlinear control affine system~\eqref{eq:SysWithState},
    availability of the output signal~$y$ is not required on the whole interval~$[t_0,\hat{t}]$
    during Step~\ref{agostep:FMPCFirst} of~\Cref{Algo:DiscreteFMPC},
    but measurement of the state~$\hat{x}=x(\hat{t};t_0,x^0,u_{\rm MPC})$ is sufficient. 
\end{remark}

\begin{remark}
Note that while the time shift $\delta>0$ is an upper bound for the step length~$\tau>0$ of the  control signals,
$\delta$ is allowed to be larger than $\tau$ under the condition that the partition~$\cP$ contains $(t_0+k\delta)_{k\in\N_0}$ as a subsequence.
In this case,  several control signals are applied to the system between two steps of the MPC~\Cref{Algo:DiscreteFMPC}.
This can also be interpreted as a multi-step MPC scheme, cf.~\cite{worthmann2014role}.
\end{remark}

In the following main result we show 
that the sampled-data MPC \Cref{Algo:DiscreteFMPC}
is initially and recursively feasible for every prediction horizon~${T>0}$.
\begin{theorem}\label{Th:FunnelMPC}
Consider system~\eqref{eq:Sys} with $(f,g,\oT)\in\cN^{m}$ and initial data $y^0\in\con([t_0-\sigma,t_0],\R^m)$.
Let $\psi\in\cG$ and $\ y_{\rf}\in W^{1,\infty}(\Rp,\R^{m})$.
Then, there exists $u_{\max}>0$ and~$\tau>0$ such that 
\Cref{Algo:DiscreteFMPC} with $T>0$ and $\delta>0$ is initially and
recursively feasible, i.e., at time $\widehat t = t^0$ and at each
successor time $\widehat t\in t^0+\delta\N$ the OCP~\eqref{eq:DiscreteFMPCOCP}
has a solution.
In particular, the closed-loop system consisting of~\eqref{eq:Sys} and 
feedback~\eqref{eq:FMPC-fb}
has a (not necessarily unique) global solution $y:[t^0,\infty)\to\R^m$ and the corresponding input is given by
\[
        u_{\rm MPC}(t) = \mu(t,y|_{[t_0-\sigma,\hat t]}).
\]
Furthermore, each global solution~$x$ with corresponding input $u_{\rm MPC}$ satisfies:
\begin{enumerate}[(i)]
    \item\label{th:item:BoundedInput}
        $\fa t\ge t^0:\quad \Norm{u_{\rm MPC}(t)}\leq u_{\rm max}$.
    \item\label{th:item:ErrorInFunnel} The error $e=y-y_{\rf}$ evolves within the error boundaries
        $\cF_{\psi}$, i.e., $\Norm{e(t)} \le \psi(t)$ for all $t\ge t^0$.
\end{enumerate}
\end{theorem}

The proof is relegated to the appendix. Here, we sketch the main ideas.
First, using the construction of the sampled-data controller~\cite{LanzaDenn23sampled}, we show that there exist sampled-data inputs
satisfying the input constraints, and achieving the output constraints.
Invoking the definition of admissible cost functions~\eqref{eq:stageCostFunnelMPC}, we establish that the cost function is finite if and only if the corresponding input belongs to the set of admissible controls
\[
    \cU_{I}(u_{\max},\hat{y})\!:=\!\setdef{\!\!\!u\in L^\infty(I,\R^m)\!\!\!}
    {\!\!\!\!\!
        \begin{array}{l}
            y(t;\hat{t},\hat{y},u)\in\cD_t \!\!\fa t\in I,\\
            \SNorm{u}\le u_{\max}
        \end{array} \!\!\!\!\!\!\!
    }\!,
\]
where $\hat{y}\in\con([t_0-\sigma, \hat{t}],\R^m)$ with $\hat{y}(t)\in\cD_{t}$ for all {$t\in[t_0,\hat{t}]$}.
Finally, we show that minimal costs can be achieved by sampled-data 
controls, and in fact, that these controls 
are contained in the set of admissible controls. 

Although \Cref{Th:FunnelMPC} is an existence results, we emphasize that feasible choices for the step length~$\tau > 0$
and the maximal control~$u_{\rm max}$ are given explicitly in~\eqref{eq:umax} and~\eqref{eq:tau}.
The bound~$u_{\rm max}$ is constructed using worst-case estimates of the system dynamics and the funnel function~$\psi$.
It is large enough such that a constant control input of this magnitude can steer the system output~$y$ away the boundary of the funnel~$\cF_{\psi}$.
The step length~$\tau > 0$, on the other hand, is small enough to avoid overshooting within one sampling interval.
Since the derivation of both quantities requires some additional notation, we relegate the explicit expressions to the appendix. 

\section{Simulation} \label{Sec:Simulation}
\begin{figure}[ht]
         \centering
          \def\svgwidth{70pt}    
\begingroup%
  \makeatletter%
  \providecommand\color[2][]{%
    \errmessage{(Inkscape) Color is used for the text in Inkscape, but the package 'color.sty' is not loaded}%
    \renewcommand\color[2][]{}%
  }%
  \providecommand\transparent[1]{%
    \errmessage{(Inkscape) Transparency is used (non-zero) for the text in Inkscape, but the package 'transparent.sty' is not loaded}%
    \renewcommand\transparent[1]{}%
  }%
  \providecommand\rotatebox[2]{#2}%
  \newcommand*\fsize{\dimexpr\f@size pt\relax}%
  \newcommand*\lineheight[1]{\fontsize{\fsize}{#1\fsize}\selectfont}%
  \ifx\svgwidth\undefined%
    \setlength{\unitlength}{595.27559055bp}%
    \ifx\svgscale\undefined%
      \relax%
    \else%
      \setlength{\unitlength}{\unitlength * \real{\svgscale}}%
    \fi%
  \else%
    \setlength{\unitlength}{\svgwidth}%
  \fi%
  \global\let\svgwidth\undefined%
  \global\let\svgscale\undefined%
  \makeatother%
  \begin{picture}(1,1.41428571)%
    \lineheight{1}%
    \setlength\tabcolsep{0pt}%
    \put(0,0){\includegraphics[width=\unitlength,page=1]{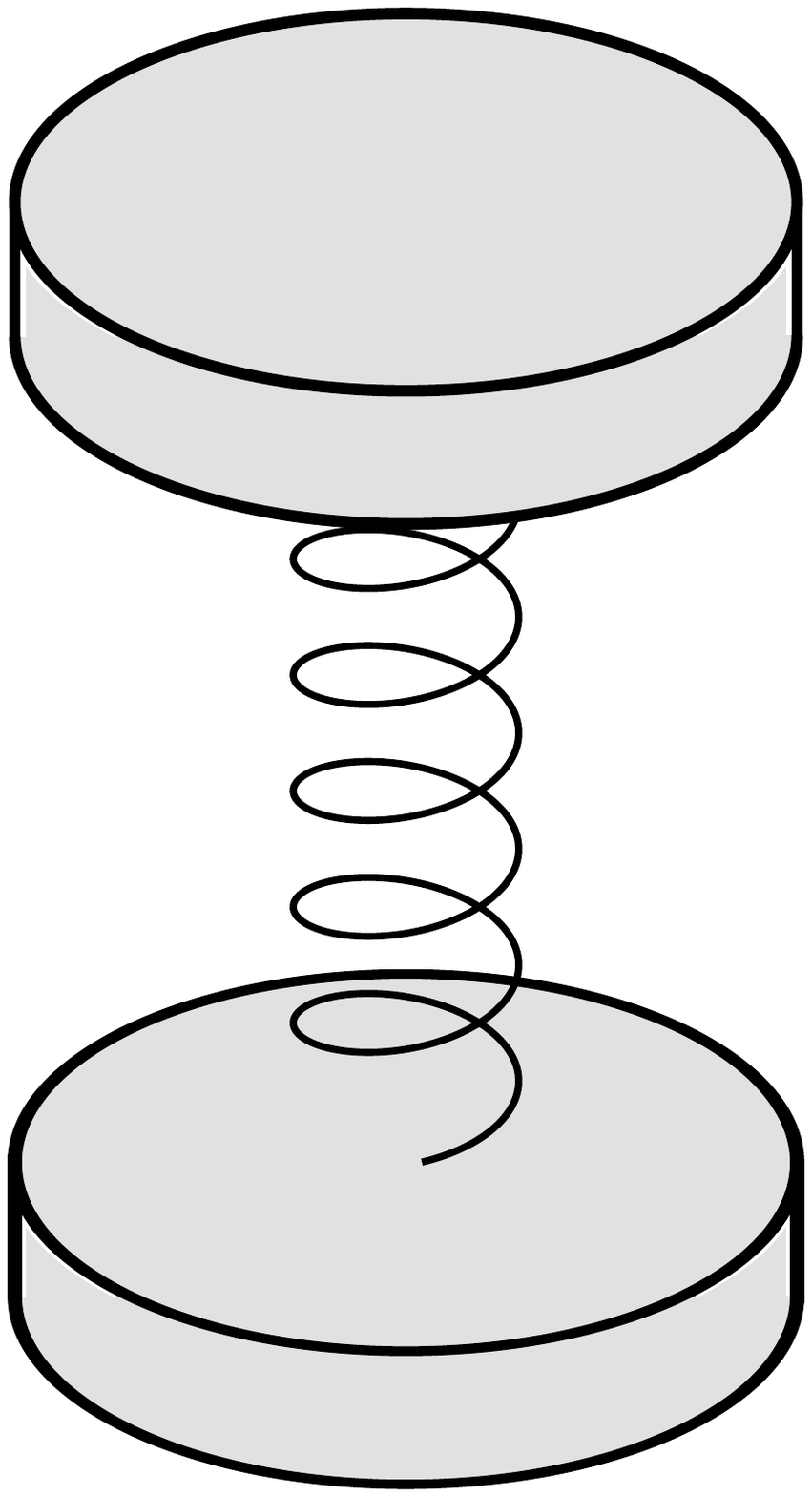}}%
    \put(0.67609275,0.49197402){\color[rgb]{0.25490196,0.46666667,0.24705882}\makebox(0,0)[lt]{\lineheight{0}\smash{\begin{tabular}[t]{l}$z_1$\end{tabular}}}}%
    \put(0.66492404,1.33212705){\color[rgb]{0.25490196,0.46666667,0.24705882}\makebox(0,0)[lt]{\lineheight{0}\smash{\begin{tabular}[t]{l}$z_2$\end{tabular}}}}%
    \put(0,0){\includegraphics[width=\unitlength,page=2]{MassChain_model.pdf}}%
    \put(0.86068449,0.27293541){\color[rgb]{0.63529412,0.10196078,0.10980392}\makebox(0,0)[lt]{\lineheight{0}\smash{\begin{tabular}[t]{l}$u$\end{tabular}}}}%
    \put(0,0){\includegraphics[width=\unitlength,page=3]{MassChain_model.pdf}}%
    \put(0.64236566,0.66745031){\color[rgb]{0,0,0}\makebox(0,0)[lt]{\lineheight{1.25}\smash{\begin{tabular}[t]{l}$k,d$\end{tabular}}}}%
  \end{picture}%
\endgroup%

         \caption{Torsional oscillator. The figure is based on \cite[Fig.~2.7]{Druecker22}, edited to the case of two flywheels for the present purpose.}
         \label{Fig:TorsOscill}
\end{figure}

The theoretical results are illustrated by a numerical example.
We consider a torsional oscillator with two flywheels, which are connected by a rod, see \Cref{Fig:TorsOscill}. %
Such a system can be interpreted as a simple model of a driving train, cf.~\cite{Pham19,Druecker22}.
The equations of motion for the torsional oscillator are given by
\begin{equation*} 
\begin{aligned}
    \begin{bmatrix} I_1 & 0 \\ 0 & I_2 \end{bmatrix}
    \begin{pmatrix} 
    \ddot z_1(t) \\ \ddot z_2(t)
    \end{pmatrix} = &
    \begin{bmatrix}
        -d & d \\ d & - d
    \end{bmatrix}
    \begin{pmatrix}
        \dot z_1(t) \\ \dot z_2(t)
            \end{pmatrix}
        \\ + &
        \begin{bmatrix}
            -k & k \\ k & -k
        \end{bmatrix}
        \begin{pmatrix}
            z_1(t) \\ z_2(t)
        \end{pmatrix}
        + \begin{bmatrix}
            1 \\ 0
        \end{bmatrix} u(t),
\end{aligned}
\end{equation*}
where for $i=1,2$ (the index~$1$ refers to the lower flywheel) $z_i$ is the rotational position of the flywheel, $I_i > 0$ is the inertia, $d,k > 0$ are damping and torsional-spring constant, respectively.
We aim to control the oscillator such that the lower flywheel follows a given velocity profile. Hence, we choose 
$
    y(t) = \dot z_1(t)
$
as the output.
To remove the rigid body motion from the dynamics, we introduce $\hat z:= z_1 - z_2$.
With this new variable, setting $x:=(\hat z, \dot z_1, \dot z_2)$ the dynamics can be written as
\begin{equation*}
    \dot x(t) = A x(t) + B u(t), \quad y(t) = C x(t)  = \dot z_1(t),
\end{equation*}
where
\begin{equation*} %
\begin{aligned}
   & M := \begin{bmatrix} 1&0&0\\ 0&I_1&0 \\ 0&0&I_2 \end{bmatrix}, \
    \tilde A := \begin{bmatrix} 0&1&-1\\-k&-d&d\\k&d&-d \end{bmatrix}, \
    \tilde B: = \begin{bmatrix}0\\1\\0\end{bmatrix}, \\
    & A : = M^{-1} \tilde A, \
    B := M^{-1} \tilde B, \ C= \begin{bmatrix} 0 & 1 & 0\end{bmatrix}.
    \end{aligned}
\end{equation*}
Using standard techniques, see, e.g.,~\cite{ilchmann1991non}, and invoking \Cref{Rem:StateRepresentation}, the reduced dynamics of the torsional oscillator can then be written in input/output form
\begin{equation} \label{eq:TorsOscill-IO}
    \begin{aligned}
        \dot y(t) &= R y(t) + S \eta(t) + \Gamma u(t), \\
        \dot \eta(t) &= Q \eta(t) + P y(t),
    \end{aligned}
\end{equation}
where~$\eta$ is the internal state, and 
$R = \frac{-d}{I_1}, \ S = \frac{1}{I_1} \begin{smallbmatrix} k & d \end{smallbmatrix}$,  
$Q = \frac{1}{I_2} \begin{smallbmatrix} 0 & I_2 \\ -k & -d \end{smallbmatrix}$,
$P = \frac{1}{I_2} \begin{smallbmatrix} -I_2 \\d \end{smallbmatrix}$.
Note that~$Q$ is a stable matrix, i.e., its eigenvalues are on the left half plane. Thus, the internal dynamics are bounded-input bounded-state stable.

\begin{figure}
    \centering
    \includegraphics[scale=0.4]{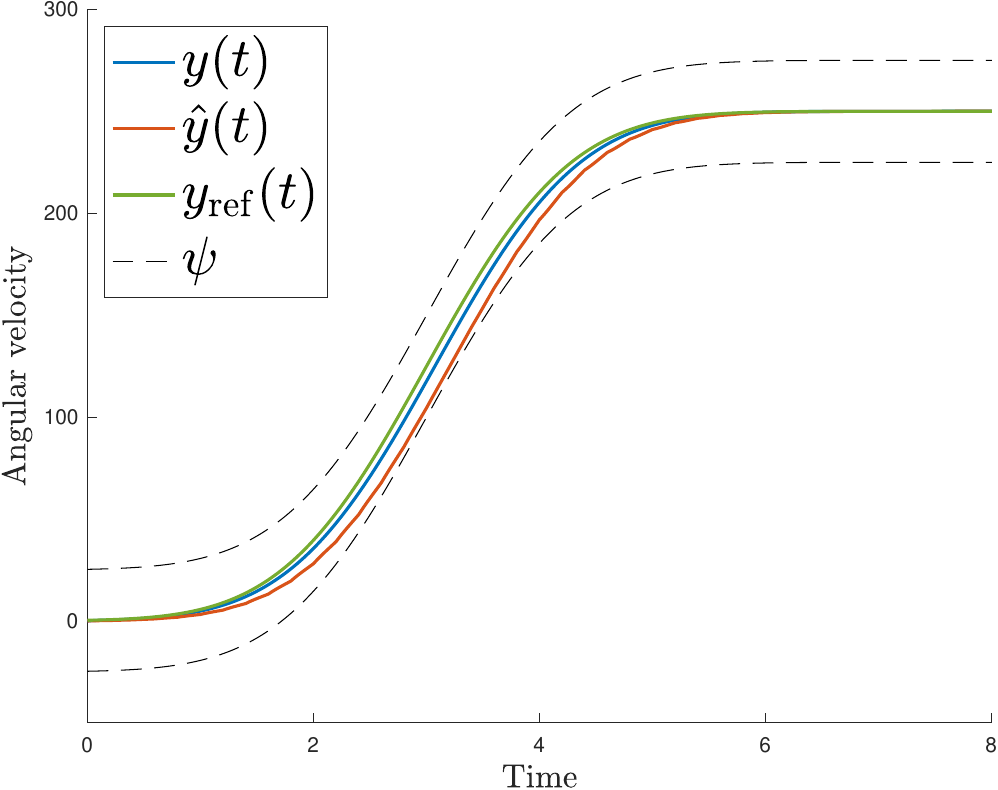}
    \caption{Outputs and reference, with error boundary.}
    \label{Fig:Outputs}
\end{figure}

\begin{figure}
    \centering
    \includegraphics[scale=0.4]{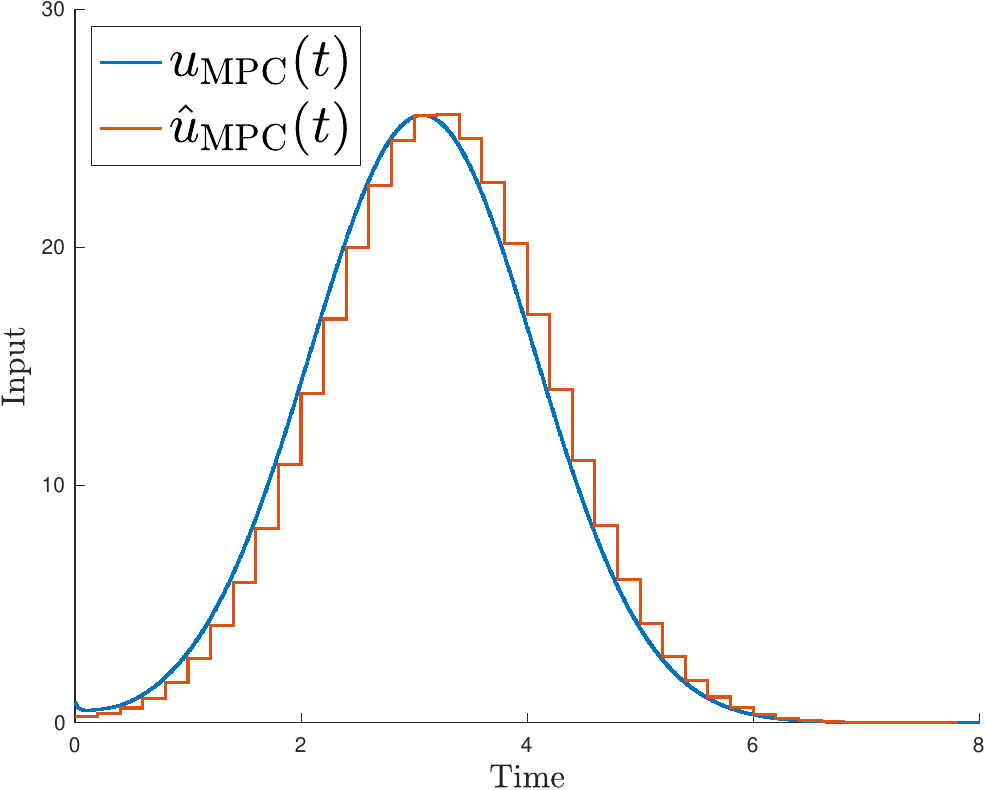}
    \caption{Control inputs.}
    \label{Fig:Control}
\end{figure}

The high-gain matrix is given by $\Gamma := CB = 1/I_1 > 0$.
For purpose of simulation, we choose the reference
\begin{equation*}
    y_{\rm ref}(t) = \frac{250}{2} \left(1 + \frac{1}{\sqrt{2\pi}} \int_0^t \me^{-\frac{1}{2}(s-3)^2} \text{d}s \right),
\end{equation*}
which is a modified version of the error function ($\textsc{erf}$) and represents a smooth transition from zero rotation to a (approximately) constant angular velocity of~$250$ rotations per unit time. Thus, $\|y_{\rm ref}\|_\infty \le 250$, 
$\|\dot y_{\rm ref}\|_\infty = 250/\sqrt{2\pi}$.
Inserting the dimensionless parameters
$I_1 = 0.136$, $I_2 = 0.12$, $k=10$, and $d=16$, and invoking the reference~$y_{\rm ref}$ and the constant error tolerance~$\psi = 25$ (we allow~$10\%$ deviation), we may derive worst case bounds on the system dynamics by estimating the explicit solution of the linear equations~\eqref{eq:TorsOscill-IO}. 
We compute these bounds in order to estimate a sufficiently large $u_{\max}>0$ as in~\eqref{eq:umax}.
For the sake of simplicity, we will assume $\eta(0) = 0$, which does not cause loss of generality.
For $\|y\|_\infty \le \|y_{\rm ref}\|_\infty + \psi$ we estimate
\begin{equation*}
\begin{aligned}
  \forall \, t \ge 0\, :  \ \   \| \eta(t) \| %
    & \le \frac{M}{\mu}\|P\| ( \|y_{\rm ref}\|_\infty + \psi) ;%
\end{aligned}
\end{equation*}
where $M := \sqrt{\|K^{-1}\|\|K\|}$ and $\mu := 1/(2\|K\|)$, and $K \in \R^{2 \times 2}$ solves the Lyapunov equation $KQ + Q^\top K + I(2) = 0$, where $I(2)$ is the two dimensional identity matrix.
Inserting the values, we find that the estimates for step length of the control signal~$\tau$
and maximal control provided in~\eqref{eq:umax},~\eqref{eq:tau} are satisfied with $\tau = 0.0047$, and $u_{\rm max} = 358$.
We choose the time shift $\delta=\tau$, i.e., a constant control is applied to the system between two iterations in \Cref{Algo:DiscreteFMPC}.
Further, the prediction horizon is set as~$T=10 \delta$.
For the purpose of simulation, we use the cost function~\eqref{eq:QuadraticCost} with~$\lambda_u = 10^{-1}.$
The results are depicted in \Cref{Fig:Outputs,Fig:Control}.
We stress that the estimates in~\eqref{eq:umax},~\eqref{eq:tau} are very conservative.
To demonstrate this aspect, we run a second simulation, where we chose $\tau = 0.2$, and $T=1$.
The results of this simulation are labeled as~$\hat y, \hat u$, respectively.
With this much larger uniform step length, the tracking objective can be satisfied as well, cf. \Cref{Fig:Outputs,Fig:Control}.
Note that the maximal applied control value is much smaller than the (conservative) estimate~$u_{\rm max}$ which satisfies~\eqref{eq:umax}.
The simulations have been performed with \textsc{Matlab} using the \textsc{CasADi}\footnote{http://casadi.org} framework~\cite{Andersson2019}.

\section{Conclusion} \label{Sec:Conclusion}

We proposed an MPC scheme achieving output tracking within prescribed bounds extending the work~\cite{BergDenn21} to sampled-data systems.
Moreover, we provided explicit, so far, conservative bounds on the required sampling frequency and control effort.
It is still an open question how this conservatism can be reduced and how a priori given bounds on the control input 
can be incorporated in the MPC scheme while preserving inter-sampling tracking guarantees.

Moreover, future research may extend this research to systems with arbitrary relative degree or towards a robustification similar to~\cite{BergDenn23a}, again invoking the key results~\cite{LanzaDenn23sampled}.

\appendix \label{Sec:Proof}
Throughout the appendix, let the assumptions of~\Cref{Th:FunnelMPC} hold.
In preparation of proving \Cref{Th:FunnelMPC}, we note some observations for later use
and recall some results from~\cite{LanzaDenn23sampled} adapted to the current setting.
To achieve that the tracking error~$e:= y-y_{\rf}$ evolves within~$\cF_\psi$ with $\psi\in\cG$, 
it is necessary that the output~$y(t)$ of the system~\eqref{eq:Sys} is at every time $t\geq t_0$ an element of the set~$\cD_{t}$ as in~\eqref{eq:Def-Dt}.
For a $L^\infty$--control function $u$ bounded by $u_{\max}>0$ to achieve the control objective on an interval $I\subset\Rp$ 
with $\hat{t}:=\min I$, i.e., ensuring that the tracking error~$e$ evolves within~$\cF_\psi$,
it is necessary for $u$ to be an element of the set
\[
    \cU_{I}(u_{\max},\hat{y})\!:=\!\setdef{\!\!\!u\in L^\infty(I,\R^m)\!\!\!}
    {\!\!\!\!\!
        \begin{array}{l}
            y(t;\hat{t},\hat{y},u)\in\cD_t \!\!\fa t\in I,\\
            \SNorm{u}\le u_{\max}
        \end{array} \!\!\!\!\!\!\!
    }\!,
\]
where $\hat{y}\in\con([t_0-\sigma, \hat{t}],\R^m)$ with $\hat{y}(t)\in\cD_{t}$ for all $t\in[t_0,\hat{t}]$.
Consequently, for a step function with partition $\cP$ to achieve the control objective it has to be an element of
\[
    \cU_{I}^{\cP}(u_{\max},\hat{y}):= \cT_{\cP}(I,\R^m)\cap  \cU_{I}(u_{\max},\hat{y}).
\]
A solution of the system~\eqref{eq:Sys} which fulfills the control objective up to a time $\nu>t_0$ is an element of the set 
\[
    \cY_\nu\!:=\!\setdef
    {\!\!\zeta\!\in\! \con([t_0-\sigma,\infty),\R^m)\!\!\!}
        {\!\!\!\!
        \begin{array}{l}
             \zeta =y^0, \! \fa t\!\in\! [t_0,\nu]\!:\\  \zeta(t)\!-\!y_{\rf}(t)\!\in\!\cD_t
        \end{array}\!\!\!\!\!
        }\!.
\]
This is the set of all functions $\zeta\in\con([t_0-\sigma,\infty),\R^m)$ which coincide with $y^0$
on $[t_0-\sigma,t_0]$ and evolve within the funnel $\cD_{t}$ on the interval $[t_0,\nu]$.
With this notation, we may state the following existence result, which is a particular version of~\cite[Lem.~2.2]{LanzaDenn23sampled}.
\begin{lemma} \label{Lemma:DynamicsBounded}
Under the assumptions of~\Cref{Th:FunnelMPC}, there exist constants $f_{\max}$, $g_{\max}$, $g_{\min}>0$ such that for every 
$\nu> t_0$, $\zeta\in\cY_{\nu}$, $z\in\R^{m}\backslash\cbl0\cbr$ and $t\in [t_0,\nu)$ 
\begin{equation} \label{eq:fmax_gmax_gmin}
\begin{aligned}
    f_{\max} &\geq \SNorm{f((\oT(\zeta))|_{[t_0,\nu)})} ,\\
    g_{\max} &\geq \SNorm{g((\oT(\zeta))|_{[t_0,\nu)})} ,\\
    g_{\min} &\leq \frac{\al z, g((\oT(\zeta))|_{[t_0,\nu)}(t))z\ar }{\Norm{z}^2}.
\end{aligned}
\end{equation}
\end{lemma}
Note that the existence of $g_{\min}$ is a direct consequence of the positive definiteness of the function $g$ as assumed in~\eqref{Ass:g_sign_definite}.
In virtue of \Cref{Lemma:DynamicsBounded} let
\begin{equation*} 
\begin{aligned}
\kappa_0 &:= \SNorm{\psi \dd{t}1/\psi} + \SNorm{1/\psi} ( f_{\max} + \|\dot{y}_{\rf}\|_\infty ),\\
\beta &>\tfrac{2\kappa_0}{g_{\min}\SNorm{\psi}},\\
\kappa_1 &:= \kappa_0+\SNorm{1/\psi}g_{\max}\beta, 
\end{aligned}
\end{equation*}
and $\tau\in(0,\kappa_0/\kappa_1^2]$.
It has been proven in~\cite[Thm.~3.1]{LanzaDenn23sampled} that 
the application of the ZoH control
\begin{equation} \label{eq:ZoHControl}
     \, u_{\rm ZoH}(t) =
    \begin{cases}
        0, \!\! &  \|  e(t_i)\| < \psi(t_i) \left( 1-\tfrac{\kappa_0^2}{\kappa_1^2} \right), \\
         - \beta  \tfrac{ \psi(t_i) e(t_i)}{\|e(t_i)\|^2}  , \!\! &\|  e(t_i)\| \ge \psi(t_i) \left( 1-\tfrac{\kappa_0^2}{\kappa_1^2} \right) ,
    \end{cases}
\end{equation}
to the system~\eqref{eq:Sys}, for $t \in  [t_i, t_i + \tau)$ with $t_i=i\tau$, $i\in\N_0$, 
yields $y(t;t_0,y^0,u_{\rm ZoH})\in\cD_{t}$ for all $t\in[t_0,\infty)$.
The key aspect here is that~$\beta > 0$ is large enough to compensate worst-case dynamical behavior, and~$\tau > 0$ is small enough (but fixed) to avoid overshooting in one sampling interval.
Note that $\SNorm{u_{\rm ZoH}}\leq\tfrac{\beta}{1-\kappa_0^2/\kappa_1^2}$.
Consequently 
    $\cU_{[t_0,\infty)}^{\cP}(u_{\max},y_0)\neq \emptyset$
for $u_{\max} \ge \tfrac{\beta}{1-\kappa_0^2/\kappa_1^2}$
and all partitions $\cP$ of $[t_0,\infty)$ with $\Abs{P}\leq\kappa_0/\kappa_1^2$.
If, for $\hat{t}\in \tau\N$, a on the interval $[t_0,\hat{t}]$ piece-wise constant control $u\in\cU_{[t_0,\hat{t}]}^{\cP}(u_{\max},y^0)$ is applied
to the system~\eqref{eq:Sys} on the interval $[t_0,\hat{t}]$, then the tracking error~$e$
evolves within the funnel~$\cF_{\psi}$, i.e., $y(t;t_0,y^0,u)\in\cD_{t}$
for all $t\in [t_0,\hat{t}]$. In particular, $y(\hat{t};t_0,y^0,u)\in\cD_{\hat{t}}$.
Then, $y(\cdot;t_0,y^0,u)$ can be extended to a function $\zeta\in\con([t_0-\sigma,\infty),\R^m)$, i.e., 
$\zeta|_{[t_0,\hat{t}]}\equiv y(\cdot;t_0,y^0,u)$. The function $\zeta$ is an element of $\cY_{\hat{t}}$.
Therefore, the prerequisites of \cite[Thm.~3.1]{LanzaDenn23sampled} are fulfilled with the same constants as in \eqref{eq:ZoHControl}.
This means that the application of $u_{\rm ZoH}$ to the system~\eqref{eq:Sys} on the interval 
$[\hat{t},\infty)$ ensures that the tracking error~$e$ evolves within the funnel~$\cF_{\psi}$ 
on the interval $[\hat{t},\infty)$. 
Thus, we have $u_{\rm ZoH}\in\cU_{[t_0,\hat{t}]}^{\cP}(u_{\max},y(\hat{t};t_0,y^0,u))$.

Summing up our observations, we state the following direct consequence of~\cite[Lem.~2.2, 
Thm.~3.1]{LanzaDenn23sampled}.
\begin{corollary}\label{Corollary:RecursiveControl}
    Under the assumptions of~\Cref{Th:FunnelMPC}, let 
    \begin{equation} \label{eq:umax}
     u_{\max} \ge\tfrac{\beta}{1-\kappa_0^2/\kappa_1^2}   
    \end{equation}
     and $\cP=(t_k)_{k\in\N_0}$
    be a partition of the interval $[t_0,\infty)$ with 
    \begin{equation} \label{eq:tau}
    \Abs{\cP}\leq\kappa_0/\kappa_1^2  .
    \end{equation}
    Then, $\cU_{[t_0,\infty)}^{\cP}(u_{\max},y_0)\neq \emptyset$.
    Furthermore, for all $k,j\in\N_0,k>j$, we have
    \begin{align}\label{eq:ExistenceRecurisveControls}
        \fa u\in\cU_{[t_0,t_j]}^{\cP}(u_{\max},y^0)\!:
        \cU_{[t_0,t_k]}^{\cP}(u_{\max},y(t_j;t_0,y^0,u))\neq\emptyset.
    \end{align}
\end{corollary}
With these preliminaries at hand, we may now prove~\Cref{Th:FunnelMPC}. 

\noindent
\emph{Proof of~\Cref{Th:FunnelMPC}}:
According to \Cref{Corollary:RecursiveControl}, there exist~$u_{\max}>0$ and~$\tau>0$ such that  
$\cU_{[t_0,\infty)}^{\cP}(u_{\max},y^0)\neq\emptyset$ for a partition~$\cP$ with $\Abs{\cP}\leq\tau$.
Let $T\geq\delta$ be arbitrary but fixed.
$\cU_{[t_0,\infty)}^{\cP}(u_{\max},y^0)\neq\emptyset$ implies that 
 $\cU_{[t_0,t_0+T]}^{\cP}(u_{\max},y^0)$ is non-empty as well.
Fact~\eqref{eq:ExistenceRecurisveControls} implies that if, at every time instance~$\hat{t}\in t_0+\delta\N_0$
in Step~\ref{algo:item:BoundedInput} of~\Cref{Algo:DiscreteFMPC}, a control $u\in\cU_{[\hat{t},\hat{t}+T]}^\cP(u_{\max},\hat{y})$
is applied to the system~\eqref{eq:Sys},
where $\hat{y}=y(\hat{t};t_0,y^0,u_{\rm MPC})$ is the output of the system at time~$\hat{t}$, 
then $\Norm{e(t)}=\Norm{y(t)-y_{\rf}(t)}\leq\psi(t)$ for all $t\in[\hat{t},\hat{t}+\delta]$ and 
$\cU_{[\hat{t}+\delta,\hat{t}+\delta+T]}^\cP(u_{\max},y(\hat{t}+\delta;\hat{t},\hat{y},u))\neq\emptyset$.
In particular, this inductively implies that claims (\ref{th:item:BoundedInput}) and 
(\ref{th:item:ErrorInFunnel}) of~\Cref{Th:FunnelMPC} are fulfilled.

Therefore, it only remains to show that if~$\cU_{[\hat{t},\hat{t}+T]}^\cP(u_{\max},\hat{y})$ is non-empty for
some $\hat{t}\in t_0+\delta\N_0$ and $\hat{y}\in \con([t_0-\sigma,\hat{t}],\R^m)$ with
$\hat{y}(t)\in\cD_{t}$ for all $t\in [t_0,\hat{t}]$, then the OCP~\eqref{eq:DiscreteFMPCOCP}
has a solution~$u^{\star}\in\cU_{[\hat{t},\hat{t}+T]}^\cP(u_{\max},\hat{ y})$.
To prove this, assume $\cU_{[\hat{t},\hat{t}+T]}^\cP(u_{\max},\hat{y})\neq\emptyset$ 
for $\hat{t}\in t_0+\delta\N_0$ and $\hat{y}\in \con([t_0-\sigma,\hat{t}],\R^m)$ with
$\hat{y}(t)\in\cD_{t}$ for all $t\in[t_0,\hat{t}]$.
Define the function $J:L^{\infty}([\hat{t},\hat{t}+T],\R^m)\to \R\cup\cbl\infty\cbr$ by
\[
    J(u)=\int_{\hat{t}}^{\hat{t}+T}
    \ell(t,y(t;\hat{t},\hat{y},u), u)
    \d t.
\]
\noindent
\emph{Step 1}:
Adapting~\cite[Thm.~4.3]{BergDenn21} to the current setting we show  that $J(u)<\infty$ 
for $u\in L^{\infty}([\hat{t},\hat{t}+T],\R^m)$ 
with $\Norm{u} \le u_{\max}$ if and only if $u\in\cU_{[\hat{t},\hat{t}+T]}(u_{\max},\hat{y})$.
Given $u\in \cU_{T}(u_{\max},\hat{y})$, it follows from the definition of $\cU_{T}(u_{\max},\hat{y})$ that
$e(t):=y(t; \hat{t}, \hat{y} ,u)-y_{\rf}(t)\in\cD_t$ for all $t\in[\hat{t},\hat{t}+T]$.
Thus,
\[
    \forall\, t\in[\hat{t},\hat{t}+T]:\ \Norm{e(t)}\leq\psi(t).
\]
Therefore, $\ell(t,y(t;\hat{t},\hat{y},u), u) =\tilde{\ell}(t,y(t;\hat{t},\hat{y},u), u)$ for all $t\in[\hat{t},\hat{t}+T]$.
Since $\tilde{\ell}$ is a continuous non-negative function, there exists an $\bar{\ell}>0$ such that 
$\tilde{\ell}(t,y,u)\leq \bar{\ell}$ for all $t\in[\hat{t},\hat{t}+T]$, $u\in\cB_{u_{\max}}$ and all $y\in\cD_{t}$.
Hence,
\begin{align*}
    J(u)=   \!\int_{\hat{t}}^{\hat{t}+T}\!\!\!\tilde{\ell}(t,y(t;\hat{t},\hat{y},u),u )\d{t}
    \leq   T\bar{\ell}<    \infty.
\end{align*}
To show the opposite direction, let $u\in L^{\infty}([\hat{t},\hat{t}+T],\R^m)$ with $J(u)<\infty$.
Assume there exists $\tilde{t}\in(\hat{t},\hat{t}+T]$ with $\Norm{e(\tilde{t})}> \psi(\tilde{t})$. 
By continuity of the involved functions, there exists $\varepsilon\in(0,\tilde{t}-\hat{t})$ with 
$\Norm{e(t)}> \psi(t)$ for all $t\in(\tilde{t}-\varepsilon,\tilde{t})$.
Thus,
\begin{align*}
    J(u)&=\int_{\hat{t}}^{\hat{t}+T}
    \ell(t,y(t;\hat{t},\hat{y},u), u)
    \d t\\
    &\geq \int_{\tilde{t}-\varepsilon}^{\tilde{t}}
    \ell(t,y(t;\hat{t},\hat{y},u), u)
    \d t = \varepsilon\cdot\infty = \infty.
\end{align*}

\noindent
\emph{Step 2}:
We prove that $\min_{u\in \cU_{[\hat{t},\hat{t}+T]}^{\cP}(u_{\max},\hat{y})}J(u)$ exists.
Since $\cU_{[\hat{t},\hat{t}+T]}^{\cP}(u_{\max},\hat{y})\subset\cU_{[\hat{t},\hat{t}+T]}(u_{\max},\hat{y})$,
the set $\cU_{[\hat{t},\hat{t}+T]}(u_{\max},\hat{y})$ is non-empty by assumption.
Since ${J(u)\geq 0}$ for all $u\in \cU_{[\hat{t},\hat{t}+T]}^{\cP}$, the infimum $J^{\star}:=\inf_{u\in \cU_{[\hat{t},\hat{t}+T]}^{\cP}(u_{\max},\hat{y})}J(u)$ exists.
Let~$(u_k)_{k\in\N_0}\in\rbl\cU_{[\hat{t},\hat{t}+T]}^{\cP}(u_{\max},\hat{y})\rbr^{\N_0}$ 
be a minimizing sequence, meaning
$J(u_k)\to J^{\star}$. 
By choice in the \Cref{Algo:DiscreteFMPC}, we known that $(t_0 + \delta k)_{k\in\N_0}$ 
is a subsequence of the partition~$\cP=(t_k)_{k\in\N_0}$.
Therefore, there exists
$M,N\in\N_0$, $N\geq M$ with $t_M=\hat{t}$ and $t_n>\hat{t}+T$ for all $n>N$.
Define $u_{i,k}:= u_k(t_i)$ for $i=M,\ldots, N$. 
For every $i=M,\ldots,N$, $(u_{i,k})_{k\in\N_0}$ is a sequence in $\R^m$ 
with $\Norm{u_{i,k}}\leq u_{\max}$ for all $k\in \N$.
Thus, it has a limit point~$u_i^{\star}\in\R^m$. 
The function $u^\star$ defined by $u^\star|_{[t_i,t_{i+1})\cap [\hat{t},\hat{t}+T]}:=u^\star_i$ 
is an element of~$\cT_{\cP}([\hat{t},\hat{t}+T],\R^m)$ with $\Norm{u^\star}\leq u_{\max}$.
Up to subsequence,~$u_k$ converges uniformly to~$u^\star$.
Let~$y_k:= y(\cdot;t_0,y^0,u_k)$ be the sequence of associated responses.
By an adaption of the Steps~2,~3 of the proof of~\cite[Thm.~4.6]{BergDenn21} 
to the current setting,  we  may infer that $(y_k)$ has a subsequence (which we do not relabel)
that converges uniformly to $y^\star(\cdot;\hat{t},\hat{y},u^\star)$.
It remains to show $u^\star\in \cU_{[\hat{t},\hat{t}+T]}^{\cP}(u_{\max},\hat{y})$. This means to show $y^\star(t)\in\cD_{t}$ for all 
$t\in[\hat{t},\hat{t}+T]$. Assume there exists $\tau\in(\hat{t},\hat{t}+T]$ with $y^\star(\tau) \notin \cD_{\tau}$, i.e.,
$\Norm{y^\star(\tau)-y_{\rf}(\tau)}>\psi(\tau)$.
There exists $\varepsilon>0$ with $\Norm{y^\star(\tau)-y_{\rf}(\tau)}>\psi(\tau)+\varepsilon$.
Since the uniform convergence of $(y_k)$ towards $y^\star$ implies pointwise convergence of $(y_k)$,
there exists $K>0$ such that $\Norm{y^\star(\tau)-y_k(\tau)}<\varepsilon$ for all $k\geq K$.
Furthermore, $\Norm{y_k(\tau)-y_{\rf}(\tau)}\leq\psi(\tau)$ since $u_k\in \cU_{[\hat{t},\hat{t}+T]}^{\cP}(u_{\max},\hat{y})$ for all $k\in\N_0$.
This raises the following contradiction for $k\geq K$
\begin{align*}
    \psi(\tau)+\varepsilon&<\Norm{y^\star(\tau)-y_{\rf}(\tau)}\\
    &\leq\Norm{y^\star(\tau)-y_k(\tau)}+\Norm{y_k(\tau)-y_{\rf}(\tau)}\\
    &\leq\varepsilon+\psi(\tau).
\end{align*}
Thus, $u^\star\in \cU_{[\hat{t},\hat{t}+T]}^{\cP}(u_{\max},\hat{y})$.
It remains to show that $J^\star=J(u^\star)$ and $J(u^{\star})=\min_{u\in
\cU_{[\hat{t},\hat{t}+T]}^\cP(u_{\max},\hat{y})}J(u)$, which follows along the lines of Steps~6,~7 
of the proof of~\cite[Thm.~4.6]{BergDenn21}. 
\hfill$\Box$

\section*{Acknowledgment}
We are deeply indebted to Thomas Berger (Universität Paderborn) for fruitful discussions, several helpful comments, and remarks.

\bibliographystyle{IEEEtran}
\bibliography{references}

\end{document}